\documentclass{article}
\usepackage{amsmath,amssymb,graphicx,url}
\usepackage{scalerel,stackengine,nicefrac}
\newcommand\equalhat{\mathrel{\stackon[1.5pt]{=}{\stretchto{%
    \scalerel*[\widthof{=}]{\wedge}{\rule{1ex}{3ex}}}{0.5ex}}}}

\newtheorem{theorem}{Theorem}
\usepackage[sort,nocompress]{cite}
\DeclareMathAlphabet{\mathbf}{OT1}{cmr}{bx}{it}

\begin{document}

\title{Conversions between barycentric, RKFUN, and Newton representations of rational interpolants}

\author{Steven Elsworth \and Stefan G\"{u}ttel\thanks{School of Mathematics, The University of Manchester,  M13\,9PL Manchester, United Kingdom, Emails: \texttt{\{steven.elsworth,stefan.guettel\}@manchester.ac.uk} }}
\date{}

\maketitle

\begin{abstract}
We derive explicit formulas for converting between rational interpolants in barycentric, rational Krylov (RKFUN), and Newton form. We show applications of these conversions when working with rational approximants produced by the AAA algorithm [\textsc{Y.~Nakatsukasa, O. S{\`e}te, L.\,N. Trefethen,} \emph{arXiv preprint 1612.00337,} 2016] within the Rational Krylov Toolbox and for the solution of nonlinear eigenvalue problems. 
\end{abstract}

\section{Introduction}

The Rational Krylov Toolbox (RKToolbox) is a collection of scientific computing tools based on rational Krylov techniques \cite{BeGu14}. This MATLAB toolbox implements, among other algorithms, Ruhe's rational Arnoldi method~\cite{Ruhe94} and the RKFIT method  for nonlinear rational approximation \cite{BG17b}; see also the example collection on \texttt{http://rktoolbox.org}.  At the heart of many RKToolbox algorithms are so-called RKFUNs (short for \emph{rational Krylov functions}), which are matrix pencil-based representations of rational functions \cite{BeGu15,BG17b}. The toolbox overloads more than 30~MATLAB commands for RKFUN objects using object-oriented programming, including root and pole finding (methods \texttt{roots} and \texttt{poles}), the efficient evaluation for scalar and matrix arguments (\texttt{feval}), addition, multiplication, and conversion to continued fraction form (\texttt{contfrac}). See \cite[Chapter~7]{Berl17} for more details on the RKFUN calculus.

Since version 2.7 the RKToolbox supports the conversion of rational interpolants from their barycentric form into the RKFUN format. In this brief note we explain this conversion and also establish an explicit relation to rational Newton interpolation. We then apply these connections to the problem of sampling a nonlinear eigenvalue problem, making use of the recently developed AAA algorithm \cite{nakatsukasa2016aaa}. But let us first review the various rational function representations. 

\smallskip

Assume that we are given pairs of complex numbers $(z_j, f_j)$ with the $z_j$ being pairwise distinct  ($j=0:m$). Then a rational interpolant of type $[m,m]$ can be written in \emph{barycentric form} as
\begin{equation}\label{eq:bary}
	r(z) = \frac{ \displaystyle \sum_{j=0}^m \frac{w_j f_j}{z - z_j}}{\displaystyle \sum_{j=0}^m  \frac{w_j}{z - z_j}},
\end{equation}
where the barycentric weights $w_j$ are nonzero but can otherwise be chosen freely (see, e.g., \cite{schneider1986some,lawrence2014stability,klein2012applications,nakatsukasa2016aaa}). It is easily verified that $r(z) \to f_j$ as $z \to z_j$, so $r$ is indeed a rational interpolant for the data $(z_j, f_j)$.

\smallskip

The \emph{RKFUN representation} of a type $[m,m]$ rational function is 
\[
r \equalhat (\underline{H_m},\underline{K_m},\mathbf{c}_m), 
\]
where $(\underline{H_m},\underline{K_m})$ is an $(m+1)\times m$ unreduced upper-Hessenberg pencil  (\emph{unreduced} meaning that $\underline{H_m}$ and $\underline{K_m}$ have no common zero entries on their subdiagonals) and $\mathbf{c}_m = [c_0,c_1,\ldots,c_m]^T \in \mathbb{C}^{m+1}$ is a coefficient vector \cite{BG17b}. The matrix pencil gives rise to a sequence of rational basis functions $r_0,r_1,\ldots,r_m$   satisfying the \emph{rational Arnoldi decomposition} \cite{Ruhe94,BeGu15}
\begin{equation}\label{eq:rad}
	z\, [ r_0(z),r_1(z), \ldots,r_m(z)] \underline{K_m} = [ r_0(z),r_1(z), \ldots,r_m(z)] \underline{H_m},
\end{equation}
and $r$ is defined as a linear combination of these basis functions:
\begin{equation}\label{eq:rkfun}
	r(z) = \sum_{j=0}^m c_j r_j(z). 
\end{equation}
Of course, at least one of the functions $r_j$ needs to be fixed in order for the decomposition \eqref{eq:rad} to uniquely determine all other basis functions. In RKToolbox the convention is that $r_0 \equiv 1$, but other choices are possible. 

\smallskip

Finally, a rational interpolant of type $[m,m]$ is in \emph{Newton form}
\[
	r(z) = \sum_{j=0}^m d_j b_j(z)
\]
if the  basis functions $b_j$ follow a recursion
\begin{equation}\label{eq:newton}
	b_0 \equiv 1, \quad b_{j}(z) = \frac{z - \sigma_{j-1}}{\beta_{j}(h_{j} - k_{j} z)} b_{j-1}(z), \quad  j = 1: m,
\end{equation}
with complex scalars satisfying $\beta_{j}\neq 0$, $|h_{j}|+|k_{j}|\neq 0$, and $\sigma_{j-1} \neq h_{j}/k_{j}$ for all~$j$. 
Note that the $b_j$ form a rational Newton basis with each $b_j$ having roots at the points $\sigma_{i-1}$ and poles
at $\xi_{i} := h_{i}/k_{i}$ for $i=1:j$.

\section{From barycentric to RKFUN format}

Our aim is to convert the rational function $r$ from the barycentric form \eqref{eq:bary} into the RKFUN format \eqref{eq:rkfun}. First we  write
\begin{equation}\label{eq:barynewt}
	r(z) = \sum_{j=0}^m f_j r_j(z), \quad r_j(z) = \frac{ \displaystyle  \frac{w_j}{z - z_j}}{ \displaystyle \sum_{i=0}^m  \frac{w_i}{z - z_i}}
\end{equation}
with the functions $r_j$ satisfying the recursion
\[
  w_{j-1} (z - z_j) r_j(z) = w_{j} (z - z_{j-1}) r_{j-1}(z), \quad j=1:m.
\]
Bringing terms containing $z$ to the left,
\[
  z ( w_{j-1}  r_j(z)  - w_j r_{j-1}) = w_{j-1} z_j r_j(z) - w_j z_{j-1} r_{j-1}(z), 
\]
and collecting these equations in matrix form yields a rational Arnoldi decomposition
\begin{equation}\label{eq:rad1}
   z\, [ r_0(z), r_1(z), \ldots, r_m(z) ] \underline{W_m}
   = 
   [ r_0(z), r_1(z), \ldots, r_m(z) ] Z_m\underline{W_m}
\end{equation}
with
\[
Z_m = \begin{small} \begin{bmatrix} 
   z_0 &  \\
    & z_1 \\
       & & \ddots \\
       &      &  & z_{m-1} &\\
       &      &  &         & z_{m}
   \end{bmatrix}\end{small}, \quad
  \underline{W_m} = \begin{small} \begin{bmatrix} 
   -w_1 &  \\
   w_0  & -w_2 \\
       & \ddots & \ddots \\
       &        & w_{m-2} & -w_m \\
       &        &         & w_{m-1}
   \end{bmatrix}\end{small}.
\]
Note that the $(m+1)\times m$ upper-Hessenberg pencil $(Z_m\underline{W_m},\underline{W_m})$ is indeed unreduced as we assume nonzero barycentric weights~$w_j$. We have therefore converted $r$ into the RKFUN representation
\begin{equation}\label{eq:rkfun1}
	r \equalhat (Z_m\underline{W_m},\underline{W_m},\mathbf{f}_m), \quad \mathbf{f}_m = [f_0,f_1,\ldots,f_m]^T.
\end{equation}
This form allows subsequent numerical operations on $r$ as implemented in the RKToolbox. For example, the following result from \cite[Section~5.2]{BG17b} summarizes how to compute the roots of $r$.

\begin{theorem}\label{thm:roots}
Given the RKFUN representation \eqref{eq:rkfun1} with nonzero coefficient vector~$\mathbf{f}_m$. Let $P$ be an invertible $(m+1)\times (m+1)$ matrix such that $P^{-1} \mathbf{f}_m = [\gamma,0,\ldots,0]^T$ for some nonzero $\gamma$. Then the generalized eigenvalues of the lower $m\times m$-part of the pencil $(P^{-1}Z_m\underline{W_m},P^{-1}\underline{W_m})$ correspond to the roots of~$r$ (with multiplicity).
\end{theorem}

In the \texttt{roots} method in the RKToolbox, the matrix $P$ is chosen as a Householder reflector $P = I_{m+1} - \sigma \mathbf{u}\mathbf{u}^H$. Note that Theorem~\ref{thm:roots} effectively enables us to find the $m$ roots of a rational function in barycentric form~\eqref{eq:bary} by solving an $m\times m$ generalized eigenvalue problem. The barycentric data $z_j,f_j,w_j$ appear explicitly in  the vector $\mathbf{f}_m$ and the matrices $Z_m,\underline{W_m}$. It is likely that the bidiagonal structure of $(Z_m\underline{W_m}, \underline{W_m})$ can be exploited in the solution of the generalized eigenproblem, but as $m$ is typically moderate this is not our focus here.

The problem of root finding for polynomials and rational functions in (barycentric) Lagrange form has attracted quite some research, with the polynomial case slightly more explored. A popular approach for rational root finding, explained for example in \cite{lawrence2014stability}, \cite[Section~2.3.3]{klein2012applications}, and also used in \cite{nakatsukasa2016aaa}, involves the solution of an $(m+2)\times (m+2)$ eigenvalue problem, giving rise to two spurious eigenvalues at infinity. An exception is the very general class of  CORK linearizations for nonlinear eigenvalue problems in \cite{CORK15}, which contains our formulation  and also leads to $m\times m$ eigenproblems (for scalar nonlinear eigenvalue problems, of which root finding is a special case). Theorem~\ref{thm:roots} spells out the conversion explicitly when $r$ is given in the form \eqref{eq:bary}, and it makes a direct link with the RKFUN root-finding approach in \cite[Section~5.2]{BG17b}.

\medskip

For purposes other than root finding (e.g., pole finding) it may be desirable to modify \eqref{eq:rad1} so that the first rational basis function $r_0$ is transformed into a constant. This can be achieved by utilizing the relation 
\[
	\sum_{j=0}^m r_j(z) = 1.
\]
Let $Q$ be an invertible $(m+1)\times (m+1)$ matrix with first column $[1,1,\ldots,1]^T$. (For example, $Q$ may be chosen as a multiple of a unitary matrix.) Then by inserting $Q Q^{-1}$ we transform \eqref{eq:rad1} into the decomposition
\[
   z\, [ \widehat r_0(z),  \widehat r_1(z), \ldots,  \widehat r_m(z) ] (Q^{-1} \underline{W_m})
   = 
   [  \widehat r_0(z),  \widehat r_1(z), \ldots, \widehat r_m(z) ] (Q^{-1} Z_m\underline{W_m}) 
\]
where now $\widehat r_0 \equiv 1$. The  equivalent to \eqref{eq:rkfun1} is
\begin{equation}\label{eq:rkfun2}
	r \equalhat (Q^{-1} Z_m\underline{W_m},Q^{-1} \underline{W_m}, Q^{-1} \mathbf{f}_m).
\end{equation}
By computing a QZ decomposition of the lower $m\times m$-part of the pencil $(Q^{-1} Z_m\underline{W_m},Q^{-1} \underline{W_m})$ we can restore the upper-Hessenberg structure and thereby obtain the RKFUN representation used in \cite{BeGu14}. The following theorem is an immediate consequence of \cite[Theorem~2.6]{BeGu14}. 

\begin{theorem}\label{thm:poles}
Given the RKFUN representation \eqref{eq:rkfun2}. Then the generalized eigenvalues of the lower $m\times m$-part of the pencil $(Q^{-1}Z_m\underline{W_m},Q^{-1}\underline{W_m})$ correspond to the poles of~$r$ (with multiplicity).
\end{theorem}

Again, the theorem requires an $m\times m$ eigenproblem to be solved and the involved matrices are products of structured matrices (for suitably chosen $Q$).
We have implemented the described method for converting a barycentric interpolant into RKFUN format in the utility function \texttt{util\_bary2rkfun} provided with RKToolbox version 2.7. After the conversion, all methods overloaded for RKFUNs can be called on the barycentric interpolant. We demonstrate this with a simple example.

\medskip

\noindent \textbf{Example 1:}  In an example taken from \cite[page~27]{nakatsukasa2016aaa} the AAA algorithm is used to compute a rational approximant $r(z)$ to the Riemann zeta function $\zeta(z)$ on the interval $[4-40i,4+40i]$. After conversion into an RKFUN, we evaluate the matrix function $r(A)\mathbf{b}\approx \zeta(A)\mathbf{b}$ for a shifted skew-symmetric matrix $A\in\mathbb{R}^{20\times 20}$ having eigenvalues in $[4-40i,4+40i]$ and a vector $\mathbf{b}$ of all ones. This evaluation uses the efficient rerunning algorithm described in \cite[Section~5.1]{BG17b} and requires no diagonalization of~$A$.  We then compute and display the relative error $\|\zeta(A)\mathbf{b} - r(A)\mathbf{b}\|_2/\|\zeta(A)\mathbf{b}\|_2$:
\begin{verbatim}
zeta = @(z) sum(bsxfun(@power,(1e5:-1:1)',-z)); 
[r,pol,res,zer,z,f,w,errvec] = aaa(zeta,linspace(4-40i,4+40i));
rat = util_bary2rkfun(z,f,w); % convert to RKFUN format
A = 10*gallery('tridiag',10); S = 4*speye(10);
A = [ S , A ; -A , S ]; b = ones(20,1);
f = rat(A, b); % approximates zeta(A)*b
[V,D] = eig(full(A)); ex = V*(zeta(diag(D).').'.*(V\b));
norm(ex - f)/norm(ex) 

  ans = 1.5199e-13
\end{verbatim}

\section{From barycentric to Newton form}\label{sec:Newton}

The recursion \eqref{eq:newton} for the Newton basis functions $b_j$ can be linearized as
\begin{equation}\label{eq:newt2}
	z\, [b_0(z),b_1(z),\ldots,b_m(z)] \underline{M_m} = [b_0(z),b_1(z),\ldots,b_m(z)] \underline{N_m}
\end{equation}
with matrices
\begin{align}\label{eq:MN}
\underline{M_m} &= \begin{small}\begin{bmatrix} 
   1 &  \\
\beta_1 k_1    & 1 \\
       &\ddots& \ddots \\
       &      & \beta_{m-1} k_{m-1} & 1 \\
       &      &  &         \beta_m k_m
   \end{bmatrix}\end{small}, \\ \nonumber
  \underline{N_m} &= \begin{small}\begin{bmatrix} 
   \sigma_0 &  \\
   \beta_1 h_1  & \sigma_1 \\
       & \ddots & \ddots \\
       &        & \beta_{m-1}h_{m-1} & \sigma_{m-1} \\
       &        &         & \beta_m h_m
   \end{bmatrix}\end{small}.
\end{align}
Comparing the  pencil $(\underline{M_m},\underline{N_m})$ in \eqref{eq:newt2} with the pencil $(Z_m\underline{W_m},\underline{W_m})$ in \eqref{eq:rad1} we find that they are of the same structure. In particular, starting from \eqref{eq:rad1} and defining for $j=1:m$ the quantities
\[
	\sigma_{j-1} := z_{j-1}, \quad \beta_j k_j := -\frac{w_{j-1}}{w_j}, \quad \beta_j h_j:= - \frac{z_j w_{j-1}}{w_j}, 
\]
we arrive at 
\[
   z\, [r_0(z),r_1(z),\ldots,r_m(z)] \underline{M_m} = [r_0(z),r_1(z),\ldots,r_m(z)] \underline{N_m}
\]
where $\underline{M_m},\underline{N_m}$ are given in \eqref{eq:MN}. 

What we have just shown is that the basis functions $r_j$ defined in \eqref{eq:barynewt} and associated with the barycentric interpolant can also be interpreted as Newton basis polynomials, except that the first function $r_0$ is not necessarily constant. This is a useful insight because numerical methods based on rational Newton interpolation may not need to be rewritten from scratch when switching to barycentric rational interpolation. The only change required is to convert the barycentric data $z_j,f_j,w_j$ into the Newton data $\sigma_j, \xi_j = h_j/k_j, \beta_j$ using the explicit formulas provided above. 

\section{An application to nonlinear eigenproblems}

Let $F:\Omega \to \mathbb{C}^{N\times N}$ be a matrix-valued function analytic on a domain $\Omega\subseteq \mathbb{C}$. We wish to compute points $\lambda\in\mathbb{C}$, the \emph{eigenvalues of $F$},  at which $F(\lambda)$ is singular. This is the so-called \emph{nonlinear eigenvalue problem} and many numerical methods have been developed for its solution; see, e.g., \cite{mevo04,GT17} and the references therein. In this note we will focus on the NLEIGS method described in \cite{GVMM14}, which is based on an approximate expansion of $F$ into a rational interpolant
\begin{equation}\label{eq:Rm}
	R(z) = \sum_{j=0}^m b_j(z) D_j,
\end{equation}
where the $D_j\in\mathbb{C}^{N\times N}$ are the matrix analogue to divided differences and the $b_j$ are the rational Newton basis functions defined by \eqref{eq:newton}. The rational interpolant $R(z)$ is then linearized into a (large, sparse and structured) matrix pencil $(\mathbf{A}_m,\mathbf{B}_m)$ of size $Nm\times Nm$. For completeness, we quote \cite[Theorem~3.2]{GVMM14}.

\begin{theorem}\label{thm:nleigs}
\label{th:linearization}
Given a rational eigenvalue problem \eqref{eq:Rm} with basis functions~$b_j$ defined by \eqref{eq:newton}.
Then the linear pencil
\begin{equation}\nonumber
\mathbf L_m(z) = \mathbf A_m - z \mathbf B_m
\label{eq:pencil}
\end{equation}
with
\begin{align*}
\mathbf A_m &= \begin{small}\begin{bmatrix}
h_m D_0 & h_m D_1 & \cdots & h_{m} D_{m-2} & h_{m}D_{m-1}-h_m\sigma_{m-1} D_{m}/\beta_m \\
 h_1\sigma_0 I & h_1\beta_1 I \\
 & \ddots & \ddots \\
 & & h_{m-2} \sigma_{m-3} I & h_{m-2}\beta_{m-2} I & \\
 & & & h_{m-1} \sigma_{m-2} I & h_{m-1}\beta_{m-1}I
\end{bmatrix},\end{small} 
\\
\mathbf B_m &= \begin{small}\begin{bmatrix} 
k_m D_0 & k_m D_1 & \cdots & k_m D_{m-2} & k_m D_{m-1} - h_m D_m/\beta_m \\
 h_1I & k_1\beta_1 I \\
 & \ddots & \ddots \\
 & & h_{m-2}I & k_{m-2}\beta_{m-2} I & \\
 & & & h_{m-1}I & k_{m-1}\beta_{m-1} I
\end{bmatrix},\end{small} 
\end{align*}
is a strong linearization of $R(z)$. If $(\lambda,\mathbf{x})$, $\mathbf{x}\neq \mathbf{0}$, is an eigenpair of $R$, that is $R(\lambda)\mathbf{x} = \mathbf{0}$, then $(\lambda,b(\lambda) \otimes \mathbf x)$ with $b(\lambda) := [ b_0(\lambda) , b_1(\lambda),  \ldots , b_{m-1}(\lambda)]^T$ is an eigenpair of $\mathbf L_m(z)$. Conversely, if $(\lambda,\mathbf y_m)$ is an eigenpair of $\mathbf L_m(z)$ then there exists a vector $\mathbf x$ such that $\mathbf y_m = b(\lambda) \otimes \mathbf x$ and $(\lambda,\mathbf x)$ is an eigenpair of $R$.
\end{theorem}

After linearization of $R(z)$ it  remains to compute the  eigenvalues of $(\mathbf{A}_m,\mathbf{B}_m)$ inside a   \emph{target set} $\Sigma$, a compact subset of $\Omega$ in which the desired eigenvalues of the nonlinear eigenproblem~$F$ are located. If $R\approx F$ uniformly on $\Sigma$, one can show that these eigenvalues are good approximations to the eigenvalues of $F$, and that there are \emph{no spurious eigenvalues} inside $\Sigma$ (i.e., eigenvalues of $R$ which are not eigenvalues of $F$). Depending on the size $Nm$ the linear eigenvalue problem can be solved either directly or iteratively (e.g., by a rational Krylov method \cite{Ruhe94}).

Note that the basis functions $b_j$ in \eqref{eq:newton} depend on the sampling points $\sigma_{i-1}$, the poles $\xi_i = h_i/k_j$, and the scaling parameters $\beta_i$ ($i,j=1:m$), all of which have to be chosen so that $R \approx F$ uniformly on $\Sigma$. The NLEIGS sampling procedure described in \cite[Section~5]{GVMM14} serves this purpose, using greedy Leja--Bagby points as the parameters \cite{bagby1969interpolation}. However, this approach requires the user to specify candidates for the pole parameters $\xi_i$, which may be a disadvantage in particular when the analytic form of $F$ is not accessible. 

In order to overcome this drawback, we consider instead of \eqref{eq:Rm} the barycentric form 
\begin{equation}\label{eq:baryR}
	R(z) = \frac{ \displaystyle \sum_{j=0}^m \frac{w_j}{z - z_j} F(z_j)}{\displaystyle \sum_{j=0}^m  \frac{w_j}{z - z_j}}
\end{equation}
and aim to choose the parameters $z_j,w_j$ so that $R \approx F$ uniformly on $\Sigma$. The AAA algorithm \cite{nakatsukasa2016aaa} seems to provide a powerful tool for this, except that it is only applicable to scalar functions. However, if we use instead of the matrix-valued function $F$ a \emph{scalar surrogate function}~$f$ having  the same region of analyticity and similar behaviour over $\Omega$, we expect that the sampling parameters~$z_j$ and the barycentric weights $w_j$ computed for $f$ are also good choices for interpolating the original function $F$ via~\eqref{eq:baryR}. Here we consider the scalar surrogate function $f(z) = \mathbf{u}^H F(z) \mathbf{v}$, where $\mathbf{u},\mathbf{v}\in\mathbb{C}^N$ are random vectors of unit length. 

Due to the direct link between the barycentric and Newton forms established in Section~\ref{sec:Newton}, we can modify the existing NLEIGS implementation in RKToolbox with minimal effort. Given the boundary of the target set $\Sigma$ and the functionality to evaluate $F$ in points on that boundary, the AAA--NLEIGS combination works as follows:

\begin{enumerate}
\item Discretize the boundary of the target set $\Sigma$ with sufficiently many candidate points, collected in a set $Z$.
\item Run the AAA algorithm to compute a type $[m,m]$ rational interpolant $r$ of the form \eqref{eq:bary} for the surrogate $f(z) = \mathbf{u}^H F(z) \mathbf{v}$ on the set $Z$.  (The degree $m$ is determined by the built-in stopping criterion of AAA.)
\item Convert the barycentric parameters $z_j,f_j,w_j$ into the Newton parameters $\sigma_j, \xi_j = h_j/k_j, \beta_j$ using the formulas provided in Section~\ref{sec:Newton}.
\item Build the NLEIGS linearization in Theorem~\ref{thm:nleigs} using the Newton parameters. The matrices $D_j$ are not affected by the conversion and we simply have $D_j = F(z_j) = F(\sigma_j)$.
\item Solve the linearized eigenproblem using, e.g., a rational Krylov iteration.
\end{enumerate}

We illustrate this algorithm with the help of a numerical example.

\medskip

\noindent \textbf{Example 2:} Consider the nonlinear eigenproblem $F(z) = A - z^{1/2} I$, where $A\in\mathbb{R}^{20\times 20}$ is the shifted skew-symmetric matrix from Example~1. The eigenvalues of $F$ are the squares of the eigenvalues of $A$. As target set $\Sigma$ we choose a disk of radius $50$ centered at  $10 + 50i$, inside of which are three eigenvalues of~$F$. The boundary circle is discretized by $100$ equispaced points, providing the set $Z$ of candidate points. 

\begin{figure}[ht] 
\hspace*{-.4cm}\includegraphics[width=7cm]{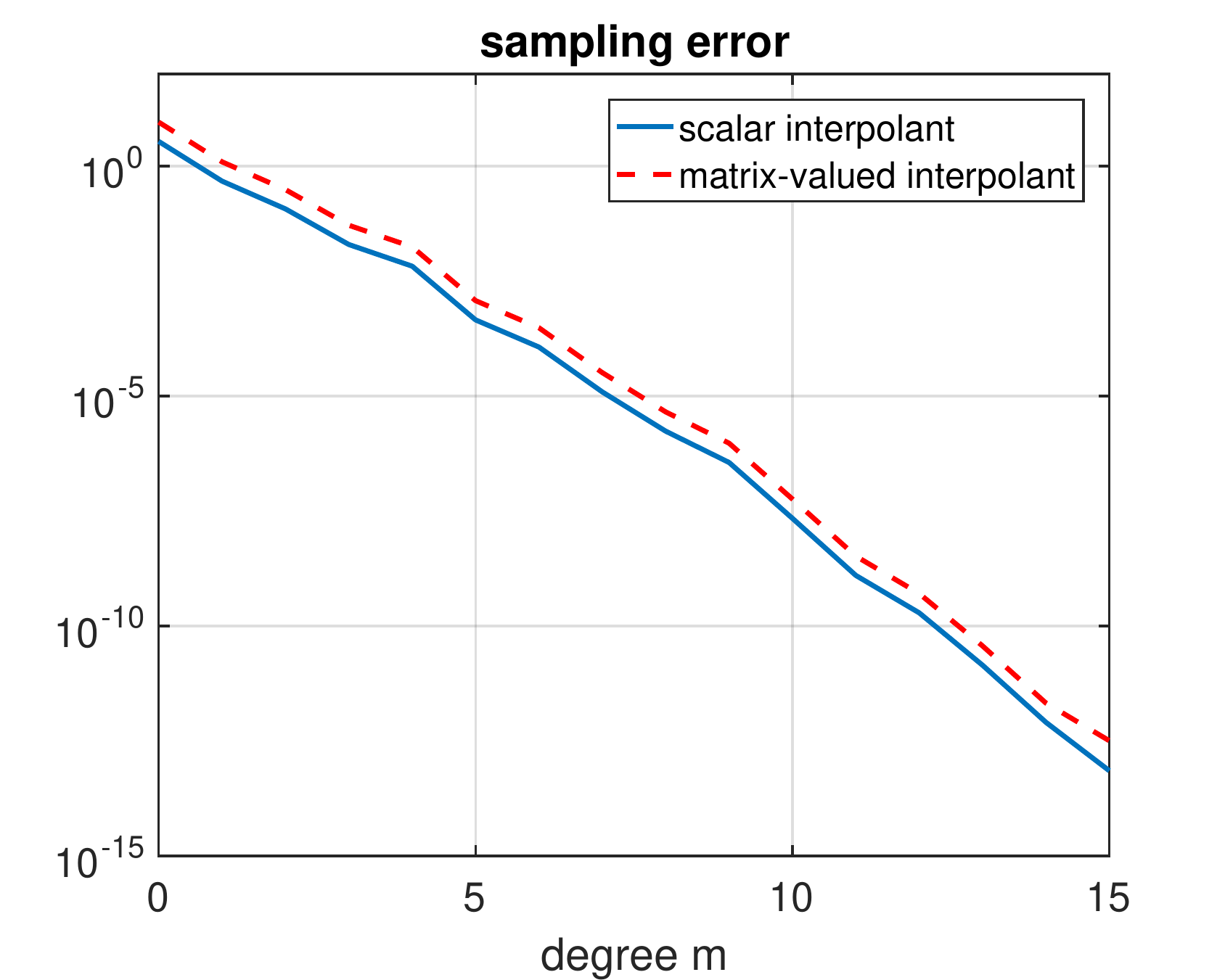}
\hspace*{-.6cm}\includegraphics[width=7cm]{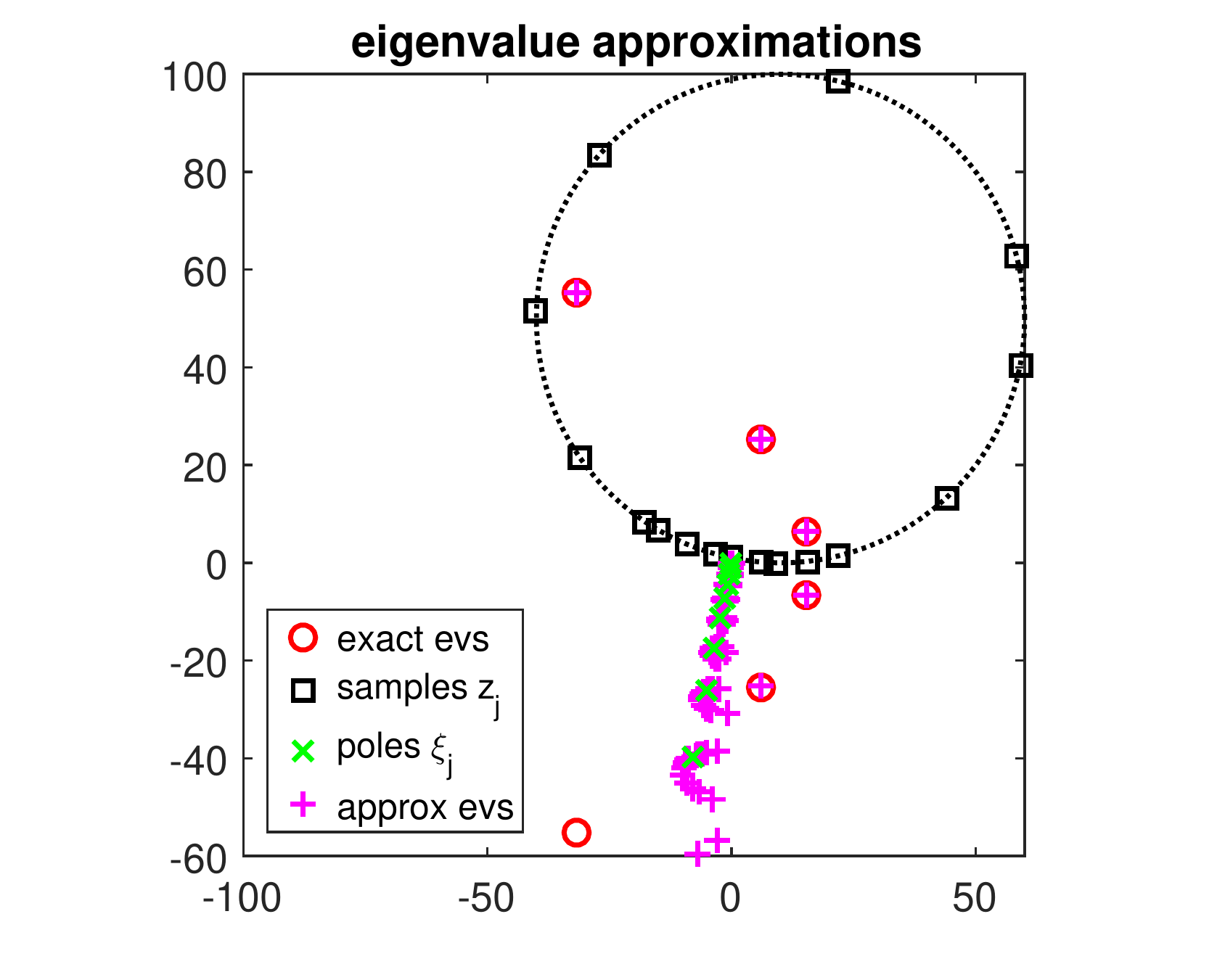} 
\caption{Solution of a nonlinear eigenvalue problem using the AAA algorithm combined with NLEIGS. Left: Accuracy of the AAA sampling procedure as the degree $m$ of the rational interpolant increases. Right: Exact eigenvalues and their approximations, together with the sampling points $z_j$ and the poles $\xi_j$ of the rational interpolant.\label{fig:nep}}
\end{figure}

We apply the AAA algorithm to $f(z) = \mathbf{u}^H F(z) \mathbf{v}$, which returns with a rational interpolant $r$ of degree $m = 15$ in barycentric form. The accuracy of $r$ measured in the uniform norm over the set $Z$ is shown on the left of Figure~\ref{fig:nep} (solid line). We also show the error $\max_{z\in Z} \|F(z) - R(z)\|_2$ for the corresponding matrix-valued interpolant $R$ (dashed line) and confirm that both approximants converge at similar rates. 

Following the algorithm outlined above, we convert the rational interpolant into Newton form and then supply the Newton parameters to the NLEIGS linearization code in RKToolbox. The resulting linear eigenvalue problem is of dimension $20m = 300$ and, for the purpose of this demonstration, small enough to be solved directly using MATLAB's \texttt{eig}. On the right of  Figure~\ref{fig:nep} we show the exact eigenvalues of $F$ (red circles) and their approximations (magenta pluses), together with the sampling points $z_j$ (black squares) and the poles of the rational interpolant (green crosses). We observe that the eigenvalues of $F$ inside the target set $\Sigma$ are well approximated by eigenvalues of the linearization, and that there are no spurious eigenvalues inside $\Sigma$.

\section{Conclusions}
We have derived explicit formulas for converting between rational interpolants in barycentric, rational Krylov (RKFUN), and Newton form, and shown two examples where this conversion is useful within the RKToolbox. In future work we plan to extend the sampling procedure to nonsquare problems, e.g., for sampling solutions of vector-valued initial and boundary value problems, and we plan to derive probabilistic error bounds for the surrogate sampling approach. 

\bigskip

\noindent \textbf{Acknowledgements.} The authors would like to thank the Engineering and Physical Sciences Research Council (EPRSC) and Sabisu for providing SE with a CASE PhD studentship. SG is grateful to Yuji Nakatsukasa, Olivier S{\`e}te, and Nick Trefethen for their hospitality during his sabbatical in 2016, and for many interesting discussions on the AAA algorithm and nonlinear eigenvalue problems which stimulated this work.

\bibliographystyle{acm}
\bibliography{strings,references}

\def\noopsort#1{}\def\l{\char32l}\def\v#1{{\accent20 #1}}
  \let\^^_=\v\def\hbk{hardback}\def\pbk{paperback}
\begin{thebibliography}{10}

\bibitem{bagby1969interpolation}
{\sc Bagby, T.}
\newblock On interpolation by rational functions.
\newblock {\em Duke Math. J. 36}, 1 (1969), 95--104.

\bibitem{Berl17}
{\sc Berljafa, M.}
\newblock {\em Rational {K}rylov Decompositions: Theory and Applications}.
\newblock PhD thesis, The University of Manchester, Manchester, UK, 2017.
\newblock Available as MIMS EPrint 2017.6 at
  \url{http://eprints.ma.man.ac.uk/2529/}.

\bibitem{BeGu14}
{\sc Berljafa, M., and G{\"u}ttel, S.}
\newblock {A} {R}ational {K}rylov {T}oolbox for {MATLAB}.
\newblock {MIMS EPrint} 2014.56, Manchester Institute for Mathematical
  Sciences, The University of Manchester, UK, 2014.
\newblock Available for download at \texttt{http://rktoolbox.org}.

\bibitem{BeGu15}
{\sc Berljafa, M., and G{\"u}ttel, S.}
\newblock Generalized rational {K}rylov decompositions with an application to
  rational approximation.
\newblock {\em SIAM J. Matrix Anal. Appl. 36}, 2 (2015), 894--916.

\bibitem{BG17b}
{\sc Berljafa, M., and G{\"u}ttel, S.}
\newblock The {RKFIT} algorithm for nonlinear rational approximation.
\newblock {\em SIAM J. Sci. Comput. 39}, 5 (2017), A2049--A2071.

\bibitem{GVMM14}
{\sc G{\"u}ttel, S., Beeumen, R.~V., Meerbergen, K., and Michiels, W.}
\newblock {NLEIGS:} {A} class of fully rational {K}rylov methods for nonlinear
  eigenvalue problems.
\newblock {\em SIAM J. Sci. Comput. 36}, 6 (2014), A2842--A2864.

\bibitem{GT17}
{\sc G{\"u}ttel, S., and Tisseur, F.}
\newblock The nonlinear eigenvalue problem.
\newblock {\em Acta Numer. 26\/} (2017), 1--94.

\bibitem{klein2012applications}
{\sc Klein, G.}
\newblock {\em Applications of linear barycentric rational interpolation}.
\newblock PhD thesis, Universit{\'e} de Fribourg, 2012.

\bibitem{lawrence2014stability}
{\sc Lawrence, P.~W., and Corless, R.~M.}
\newblock Stability of rootfinding for barycentric {L}agrange interpolants.
\newblock {\em Numer. Algorithms 65}, 3 (2014), 447--464.

\bibitem{mevo04}
{\sc Mehrmann, V., and Voss, H.}
\newblock Nonlinear eigenvalue problems: {A} challenge for modern eigenvalue
  methods.
\newblock {\em GAMM-Mitt. 27\/} (2004), 121--152.

\bibitem{nakatsukasa2016aaa}
{\sc Nakatsukasa, Y., S{\`e}te, O., and Trefethen, L.~N.}
\newblock The {AAA} algorithm for rational approximation.
\newblock {\em arXiv preprint arXiv:1612.00337\/} (2016).

\bibitem{Ruhe94}
{\sc Ruhe, A.}
\newblock Rational {K}rylov algorithms for nonsymmetric eigenvalue problems.
  {II}. {M}atrix pairs.
\newblock {\em Linear Algebra Appl. 198\/} (1994), 283--295.

\bibitem{schneider1986some}
{\sc Schneider, C., and Werner, W.}
\newblock Some new aspects of rational interpolation.
\newblock {\em Math. Comp. 47}, 175 (1986), 285--299.

\bibitem{CORK15}
{\sc {Van Beeumen}, R., Meerbergen, K., and Michiels, W.}
\newblock Compact rational {K}rylov methods for nonlinear eigenvalue problems.
\newblock {\em SIAM J. Matrix Anal. Appl. 36}, 2 (2015), 820--838.

\end{thebibliography}

\end{document}